%&plain
%
%  This is an amstex file, but using and older version of amstex
%  the relevant macros have been included, so the whole file should
%  run under plaintex instead.
%
\input oldamstex
\def\SBIMSMark#1#2#3{
 \font\SBF=cmss10 at 10 true pt
 \font\SBI=cmssi10 at 10 true pt
 \setbox0=\hbox{\SBF Stony Brook IMS Preprint \##1}
 \setbox2=\hbox to \wd0{\hfil \SBI #2}
 \setbox4=\hbox to \wd0{\hfil \SBI #3}
 \setbox6=\hbox to \wd0{\hss
             \vbox{\hsize=\wd0 \parskip=0pt \baselineskip=10 true pt
                   \copy0 \break%
                   \copy2 \break% 
                   \copy4 \break}}
 \dimen0=\ht6   \advance\dimen0 by \vsize \advance\dimen0 by 8 true pt
                \advance\dimen0 by -\pagetotal
 \dimen2=\hsize \advance\dimen2 by .25 true in
%
%   Check for publication info
%
%  \newread\jref
  \openin2=publishd.tex
  \ifeof2\setbox0=\hbox to 0pt{}
  \else 
     \setbox0=\hbox to 3.1 true in{
                \vbox to \ht6{\hsize=3 true in \parskip=0pt  \noindent  
                {\SBI Published in modified form:}\hfil\break
                \input publishd.tex 
                \vfill}}
  \fi
  \closein2
  \ht0=0pt \dp0=0pt
 \ht6=0pt \dp6=0pt
 \setbox8=\vbox to \dimen0{\vfill \hbox to \dimen2{\copy0 \hss \copy6}}
 \ht8=0pt \dp8=0pt \wd8=0pt
 \copy8
 \message{*** Stony Brook IMS Preprint #1, #2 ***}
}

\documentstyle {oldamsppt}
\magnification=1200
\tolerance=3000
\openup 6 pt
\nologo
\SBIMSMark{1997/12}{October, 1997}{}
\bigskip
\bigskip
\bigskip
\bigskip
\topmatter
\title
The Multipliers of Periodic Points in One-dimensional Dynamics 
\endtitle
\author
 Marco Martens \footnote {Institute for Mathematical Sciences, SUNY at 
Stony Brook, Stony Brook, NY 11794-3651.} 
Welington de Melo \footnote{IMPA, Estrada Dona Castorina 110, 22460-320, Rio de Janeiro, Brazil.  Partially supported by Pronex project on Dynamical 
      Systems}
\endauthor
\date{September 30 , 1997}
\endtopmatter
\bigskip
\bigskip
\bigskip
\bigskip
\bigskip
\bigskip
\centerline {\bf Abstract.} 

\bigskip

\flushpar
It will be shown that the smooth conjugacy class of an $S-$unimodal map which does not have a periodic attractor neither a Cantor attractor is determined by the multipliers of the periodic orbits. This generalizes a result by M.Shub and D.Sullivan for smooth expanding maps of the circle.
 
\newpage

\bigskip
\centerline{\bf 1. Introduction}
\bigskip

\bigskip

\flushpar
Let $f$ and $g$ be topologically equivalent smooth interval maps. That means, there is a homeomorphism $h$ such that $h\circ f=g\circ h$. We will say that $f$ and $g$ have the {\it same multipliers} if for every periodic point $x$ with period $n$ 
$$
Dg^n(h(x))=Df^n(x).
$$ 

\flushpar
Let $\Cal{U}$ be the set of unimodal maps $f:[-1,1]\to [-1,1]$ of the form
$$
f=\phi\circ q_t,
$$
where $\phi:[-1,1]\to [-1,1]$ is an orientation preserving  $C^3$ diffeomorphism with negative Schwarzian derivative and $q_t:[-1,1]\to [-1,1]$ a {\it canonical folding map}
$$
q_t(x)=-2t|x|^\alpha+2t-1
$$
for some $t\in [0,1]$ and $\alpha>1$. The number $\alpha$ is called the {\it critical exponent} of $f$.
The maps in  $\Cal{U}$ are often called $S-$unimodal maps.

\flushpar
It has been shown in [BL], [Ma] that for every map $f\in \Cal{U}$ there exists
a unique closed invariant set $A_f\subset [-1,1]$ such that
$$
\omega(x)=A_f
$$
for almost every  $x\in [-1,1]$ (in Lebesgue sense). Here $\omega(x)$ denotes the limit set of $x$. This set $A_f$ is the {\it attractor} of $f$ as defined
by Milnor [Mi].

\flushpar
There are three possibilities: the attractor is a periodic orbit, or it is a Cantor set, or it is the orbit of a periodic interval. This work will concentrate on maps whose attractor is of the third type. 
$$
\Cal{U}_0=\{f\in\Cal{U}| \text{ the topological dimension of } A_f \ne 0\}.
$$

\proclaim{Theorem} Let $f,g\in \Cal{U}_0$ be conjugated by a homeomorphism $h$. If $f$ and $g$ have the same multipliers then they belong to the same $C^2$ conjugacy class, that is, $h:A_f\to A_g$ is a $C^2$ diffeomorphism.
\endproclaim

\flushpar
This Theorem is generalization of a result obtained by M.Shub and D.Sullivan. They proved that two smooth expanding circle maps with the same multipliers are smoothly conjugated [SS]. In [L1] it has been shown that $C^2$ unimodal maps  with Fibonnaci combinatorics and with the same multipliers are $C^1$ conjugated. The proof of the Theorem will be by joining the methods in [SS] and [Ma]. 

\flushpar
For $S-$unimodal maps with critical exponent $\alpha=2$ it has been proved in [L2] that a Cantor attractor only appears for infinitely renormalizable maps  (see also [LM]). Hence for such unimodal maps the multipliers of the periodic orbits form a complete smooth invariant whenever the map does not have a periodic attractor and only finitely many renormalization. In [BKNS] examples have been constructed of non-renormalizable unimodal maps, with high critical exponent, exhibiting Cantor attractors. A characterization of the smooth conjugacy classes of such maps is not known.

\proclaim{Conjecture} Let $f,g\in\Cal{U}$  be unimodal maps which do not have periodic attractors neither Cantor attractors. If $f$ and $g$ are $C^r$, $r\ge 2$, and have the same multipliers then they are conjugated by a $C^r$ diffeomorphism.
\endproclaim

\flushpar
In the context of unimodal maps which have a quadratic-like extension, the  method presented here will prove that having the same multipliers implies that
the conjugation is even real analytic. The proof of Lemma 2.6 and 3.3 will have to be changed slightly by using the Koebe-Lemma for univalent maps.

\bigskip

\flushpar
An appendix is added in which some basic notions are defined.

\bigskip
\centerline{\bf Acknowledgements}
\bigskip

\bigskip

\flushpar
The first author would like to thank IMPA for its kind hospitality during august
1997, the period in which this rather serendipitous result was obtained.

\bigskip
\centerline{\bf 2. Markov-Maps}
\bigskip

\bigskip

\flushpar
In section 3. the conjugacy problem between unimodal maps will be reduced
to the conjugacy problem between Markov-maps. Before we define Markov-maps we will discuss the Banach space $\text{Diff}^r_+([-1,1])$, $r\ge 2$, consisting of $C^r$ orientation preserving diffeomorphisms of $[-1,1]$. The Banach space structure is obtained by the following identification. Let $\eta:\text{Diff}^r_+([-1,1])\to C^{r-2}([-1,1])$ defined by
$$
\eta_\phi=D\ln D\phi.
$$
This bijection $\eta$ is called the {\it non-linearity}. The usual Banach space structure on $C^{r-2}([-1,1])$ can be pulled back onto  $\text{Diff}^r_+([-1,1])$. In particular 
$$
|\phi|_r=\sum_{k\le r-2} ||D^k\eta_{\phi}||_0,
$$
where $||.||_0$ denotes the $C^0$ norm on the space of continuous functions.

\flushpar
Furthermore, if $f:T\to f(T)$ is an orientation preserving diffeomorphism defined on the interval $T$ then
$$
[f|T]\in\text{Diff}^r_+([-1,1])
$$
denotes the orientation preserving diffeomorphism obtained by rescaling range and domain of $f$.

\proclaim{Definition 2.1} A map $F:\bigcup M_i\to [-1,1]$ where $\{M_i|i\in \Bbb{N}\}$ is a countable collection of closed nondegenerate oriented subintervals in $[-1,1]$ with the property
$$
|[-1,1]\setminus\bigcup M_i|=0,
$$
is called a $C^r$ Markov-map, $r\ge 3$, if
\parindent=15pt
\item{1)} $F:M_i\to [-1,1]$ is an orientation preserving diffeomorphism
for every $i\in \Bbb{N}$.
\item{2)} There exists $a>0$ such that for every $i\in \Bbb{N}$ there exists
an interval $[-1,1]\supset T_i\supset M_i$ and an extension $\hat{F}:T_i\to [-1-a,1+a]$ of $F:M_i\to[-1,1]$ which is a $C^r$ diffeomorphism with negative Schwarzian derivative.
\item{3)} There exists $K<0$ such that $|[(F|{M_i})^{-1}]|_r\le K$ for every $M_i$.
\endproclaim

\flushpar
Observe, that the maps $F|M_i$ can be orientation reversing when the usual orientation on $M_i$ is used.

\flushpar
Let $H$ be a homeomorphism conjugating the $C^r$ Markov-maps $F$ and $G$. The conjugacy $H$ {\it preserves the multipliers} iff for every periodic point $x$ of $F$ with period $n$ 
$$
DG^n(H(x))=DF^n(x).
$$
Observe that in general there are many conjugating homeomorphisms between two given Markov-maps. This is the reason why the above definition differs from the 
one concerning unimodal maps.

\proclaim{Proposition 2.2} Every conjugation between Markov-maps $F$ and $G$ which preserves the multipliers is Lipschitz. In particular, such a conjugation
is absolutely continuous with respect to the Lebesgue measure. 
\endproclaim

\proclaim{Proposition 2.3} For every $C^r$ Markov-map, $r\ge 3$, there is a
$C^r$ diffeomorphism $H:[-1,1]\to [-1,1]$ such that 
$$
F_0=H^{-1}\circ F\circ H
$$
is a $C^r$ Markov-map which preserves the Lebesgue measure.  
\endproclaim

\flushpar
The proofs needs some preparation.
Let $\Sigma$ be the one-sided symbol sequence over the alphabet $\Bbb{N}$. Take $w\in \Sigma$ and consider the finite word $w_n$ consisting of only the first $n$ symbols of $w$. There is an unique interval $I_w^n\subset [-1,1]$ such that 
\parindent=15pt
\item{1)} $F^n:I_w^n\to [-1,1]$ is a diffeomorphism.
\item{2)} $F^j(I_w^n)\subset M_{w(n-j)}$ for $j=0,1,\cdots, n-1$.

\flushpar
The map $F^n: I_w^n\to [-1,1]$ is called {\it the branch of $F^n$ with 
the combinatorics $w_n$}. The collection of intervals $\{I_w^n|w\in \Sigma\}$ is denoted by $\Cal{I}^n$. For every $w\in \Sigma$ and $n\ge 0$ let $\psi_{w,n},\phi_{w,n}\in \text{Diff}^r_+([-1,1])$ be defined by
$$
\psi_{w,n}=[(F^n|I_w^n)^{-1}]
$$
and
$$
\phi_{w,n}=[(F|I_w^{n+1})^{-1}].
$$
Observe that
$$
\psi_{w,n+1}=\phi_{w,n}\circ \psi_{w,n}.
$$

\proclaim{Lemma 2.4} The collections $\Cal{I}^n$ have the following properties.
\parindent=15pt
\item{1)} $|[-1,1]\setminus \bigcup_{w} I_w^n|=0$ for each $n\ge 0$.
\item{2)} For every $I_w^n\in \Cal{I}^n$ there exists a unique interval $T_w^n$ 
and a unique $J\in \Cal{I}^{n-1}$ such that
\parindent=30pt
\item{-} $I_w^n\subset T_w^n\subset J$ and
\item{-} $\hat{F}^n:T_w^n\to [-1-a,1+a]$ is a diffeomorphism.

\parindent=15pt
\item{3)} There is a $\delta<1$ and $C>0$ such that 
$$
\frac{|I_w^n|}{|M_{w(n)}|}, \frac{|F(I_w^n)|}{|M_{w(n-1)}|}\le C\delta^n
$$
 
\flushpar
In particular, the distortion of $F^n|I_w^n$
(and $\psi_{w,n}$) is uniformly bounded (the bound $K$ is independent of 
$w\in \Sigma$ and $n\ge 0$). Moreover, $F$ is ergodic.
\endproclaim

\demo{Proof} Statement 1) and 2) follow directly from the definition of Markov-maps. Statement 3) and the distortion statement is then a direct consequence of the Koebe-Lemma (see the Appendix). To prove the ergodicity we have to show that any given  invariant set $X$ with positive Lebesgue measure has full measure. Let $D\subset [-1,1]$ be the biggest set on which each $F^n$ is defined. This set $D$ has full measure. Take a density point of $X$ with $x\in X\cap D$. Then for each $n\ge 0$ there is  $I^n_{w_n}$ with $x\in I^n_{w_n}$. Observe
$$
\frac{|[-1,1]\setminus X|}{|[-1,1]|}
\le \lim_{n\to \infty} \frac{|F^n(I_{w_n}^n\setminus X)|}{|F^n(I_{w_n}^n)|}
\le \lim_{n\to\infty} K\frac{|I_{w_n}^n\setminus X|}{|I_{w_n}^n|}. 
$$
Here we used that the distortion of $F^n|I_{w_n}^n$ is uniformly bounded by $K$.
The intervals  $I_{w_n}^n$ shrink down to the density point $x$ of $X$, the last limit has to equal zero. We proved that $X$ has full measure.
\hfill\hfill\qed $\,\,$ (Lemma 2.4)
\enddemo

\flushpar
As a consequence of Lemma 2.4 we see that the derivatives of $\psi_{w,n}$ are uniformly bounded by $K$. In particular, we can define for each $n\ge 0$ the density 
$$
\rho_n(x)=\frac12\cdot \sum_{I_w^n\in \Cal{I}^n} |D\psi_{w,n}(x)| \cdot |I_w^n|.
$$
These densities are, because of Lemma 2.4, uniformly bounded and uniformly away from zero.  Moreover, observe that $\rho_n$ is the
density of $\mu_n=F_*^n(\frac12\lambda)$, where $F_*$ be the Perron-Frobenius operator and $\lambda$ is the Lebesgue measure on $[-1,1]$. 
Let $\mu$ be a probability measure on $[-1,1]$ which is a weak limit of the sequence
$$
\frac1n\sum_{i=0}^{n-1} \mu_i.
$$
This measure $\mu$ is an invariant measure for $F$. The next Lemma will imply that $\mu$ is absolutely continuous. The ergodicity of $F$ implies that there is only one absolutely continuous invariant measure: the above sequence of measures is actually convergent.

\proclaim{Lemma 2.5} There is a continuos function $k:[0,2)\to [1,\infty)$ with $k(0)=1$ such that for every interval $T\subset [-1,1]$ and every measurable $A\subset T$
$$
\frac{1}{k(|T|)}\frac{|A|}{|T|}\le \frac{\mu(A)}{\mu(T)}\le
k(|T|)\frac{|A|}{|T|}.
$$
\endproclaim

\demo{Proof} The Koebe-Lemma (see Appendix) implies the existence of a continuous function $k:[0,2]\to [1,\infty)$ with $k(0)=1$, such that for every interval $T\subset [-1,1]$, $I_w^n\in \Cal{I}^n$ and $x,y\in I_w^n$ with $F^n(x),F^n(y)\in T$
$$
\frac{1}{k(|T|)}\le \frac{DF^n(x)}{DF^n(y)}\le k(|T|).
$$
Hence
$$
\frac{1}{k(|T|)}\le \frac{D\psi_{w,n}(x)}{D\psi_{w,n}(y)}\le k(|T|)
$$ 
for every $x,y\in T$.

\flushpar
It is enough to prove the Lemma when the set $A\subset T$ is an interval. Let $A\subset T\subset [-1,1]$ be intervals. Then there exist for every $I_w^n$
points $a_{w,n}\in A$ and $t_{w,n}\in T$ such that
$$
\frac{\mu_n(A)}{\mu_n(T)}
=\frac{\sum_{I_w^n\in \Cal{I}^n}|D\psi_{w,n}(a_{w,n})|}
         {\sum_{I_w^n\in \Cal{I}^n}|D\psi_{w,n}(t_{w,n})|} \cdot \frac{|A|}{|T|}
$$
Using the above distortion estimate we get
$$
\frac{1}{k(|T|)} \cdot \frac{|A|}{|T|}\le
\frac{\mu_n(A)}{\mu_n(T)}
\le k(|T|) \cdot \frac{|A|}{|T|}.
$$
This bound is independent of $n$. Hence it will also hold for any weak limit
$\mu$ of $\frac1n\sum_{i=0}^{n-1} \mu_i$. The Lemma is proved.
\hfill\hfill\qed $\,\,$ (Lemma 2.5)
\enddemo

\proclaim{Lemma 2.6} The sequence $\psi_{w,n}\in \text{Diff}^r_+([-1,1])$ is a Cauchy sequence for every $w\in \Sigma$. Let $\psi_w=\lim \psi_{w,n}$. The function
$\Psi:\Sigma\to \text{Diff}^r_+([-1,1])$ with
$$
\Psi(w)=\psi_w
$$
is continuous.
\endproclaim

\demo{Proof} 
Let $w\in\Sigma$ and $n\ge 0$. Remember that $\phi_{w,n}\in \text{Diff}^r_+([-1,1])$ was
constructed such that 
$$
\psi_{w,n+1}=\phi_{w,n}\circ \psi_{w,n}.
$$
Observe that
$$
\phi_{w,n}=[(F|I_w^{n+1})^{-1}]=[(F|M_{w(n+1)})^{-1}|F(I_w^{n+1})].
$$
Hence, as an immediate consequence of Lemma 2.4(3) and definition 2.1(3) we get a uniform $C>0$ and $\delta<1$ with
$$
|\phi_{w,n}|_r\le C\delta^n.
$$

\proclaim{Claim} There exists $K>0$ such that
$$
|\psi_{w,n+1}-\psi_{w,n}|_r\le K |\phi_{w,n}|_r.
$$
\endproclaim

\demo{Proof} For every $k\ge 0$ there exists polynomials 
$P^k_j:\Bbb{R}^{k+1}\to \Bbb{R}$, $j\le k$ such that for every pair $\phi,\psi\in \text{Diff}^r_+([-1,1])$ we have
$$
D^k(\eta_{\phi}(\psi(x))D\psi(x))=
\sum_{j=0}^k D^j\eta_{\phi}(\psi(x))\cdot P^k_j(D\psi(x),D^2\psi(x),\cdots, D^{k+1}\psi(x)) ,
$$
for every $x\in [-1,1]$. Using this formula we will inductively prove the Claim. By the chain rule and the definition
$$
|\psi_{w,n+1}-\psi_{w,n}|_r=|(\eta_{\phi_{w,n}}\circ \psi_{w,n}\cdot D\psi_{w,n})|_{r-2}=\sum_{j=0}^{r-2} |D^j(\eta_{\phi_{w,n}}\circ \psi_{w,n}\cdot D\psi_{w,n})|_0.
$$
First we will prove the Claim for $r=2$. Observe
$$
\aligned
|\psi_{w,n+1}-\psi_{w,n}|_2&= |\eta_{\phi_{w,n}}\circ \psi_{w,n}\cdot D\psi_{w,n}|_0\\
                           &\le K |\eta_{\phi_{w,n}}|_0.
\endaligned
$$
We used Lemma 2.4: the maps $\psi_{w,n}$ have uniformly bounded derivative. 

\flushpar
Assume the Claim is proved for $r\ge 2$ (to do the induction step we have to assume that the Markov-map is at least $C^{r+1}$). The induction assumption implies that the sequence $\psi_{w,n}$ is a Cauchy sequence: $\psi_{w,n}\to \psi_w\in \text{Diff}^r_+([-1,1])$. In particular
$$
D^{k}\psi_{w,n}\to D^k\psi_w, k\le r
$$
in the $C^0$ topology. Now use this convergence and above the expression for derivatives $D^k(\eta_{\phi}(\psi(x))D\psi(x))$ to get a bound
$$
D^{(r+1)-2}(\eta_{\phi_{w,n}}(\psi_{w,n}(x))D\psi_{w,n}(x))\le K 
\sum_{s\le (r+1)-2} |D^s\eta_{\phi_{w,n}}|_0=K|\phi_{w,n}|_{r+1}.
$$
The constant $K$ depends on the coefficients of the polynomials $P^{(r+1)-2}_j$
and the norms of the derivatives of $\psi_{w,n}$.
\hfill\hfill\qed $\,\,$ (Claim)
\enddemo

\flushpar
The main consequence of the Claim is that every sequence $\psi_{w,n}$ is a Cauchy sequence in the $C^r$ topology. Moreover, the constants obtained in the induction steps depend only on the initial distortion bound obtained from Lemma 2.4 and the exponential decay of $|\phi_{w,n}|_r$, which is also uniform. In particular, there is a uniform constant $K>0$ such that for every pair $w',w\in \Sigma$ with $w'_n=w_n$ for some $n\ge 0$ we have
$$
|\psi_{w'}-\psi_w|_r\le K \delta^n.
$$
The function $\Psi$ is continuous.

\hfill\hfill\qed $\,\,$ (Lemma 2.6)
\enddemo

\flushpar
The Borel $\sigma-$algebra on $\Sigma$ is generated by the cylinders:
$$
[w]_n=\{v\in\Sigma| v(j)=w(j) , j\le n\}.
$$
The collection of cylinders is denoted by $\Cal{C}$. Define the function
$\nu:\Cal{C}\to [0,1]$ by
$$
\nu([w]_n)=\mu(I_w^n).
$$

\proclaim{Lemma 2.7} The function $\nu$ extends to a measure on the Borel $\sigma-$algebra of $\Sigma$. 
\endproclaim

\demo{Proof} Observe $[w]_n=\bigcup_{i} [w_ni^{\infty}]_{n+1}$. To prove that $\nu$ can be
extended to a Borel measure it is enough to show that $\nu([w]_n)=\sum_{i} \nu([w_ni^\infty]_{n+1})$. The measure $\mu$ is invariant under $F$. Hence
$$
\nu([w]_n)=\mu(I_w^n)=\mu(F^{-1}(I_w^n))=\mu(\bigcup_{i} I_{w_ni^\infty}^{n+1})=
\sum_{i}\mu(I_{w_ni^\infty}^{n+1})= \sum_{i} \nu([w_ni^\infty]_{n+1}).
$$
\hfill\hfill\qed $\,\,$ (Lemma 2.7)
\enddemo

\flushpar
Because $\Psi:\Sigma\to \text{Diff}^r_+([-1,1])$ is continuous the function $\rho:[-1,1]\to (0,\infty)$ defined by
$$
\rho(x)=\int |D\psi_{w}(x)| d\nu
$$
is $C^{r-1}$. Now we get the rather peculiar

\proclaim{Lemma 2.8} The $C^{r-1}$ function $\rho$ is the density of the invariant measure $\mu$ of $F$. 
\endproclaim

\demo{Proof} The ergodicity of $F$ implies that it is enough to show that
$\rho$ is the density of some absolutely continuous invariant measure: it has to be the density of $\mu$. First observe
$$
\rho(x)=\int |D\psi_{w}(x)| d\nu=
\lim_{n\to\infty} \sum_{I_{w,n}\in \Cal{I}^n} |D\psi_{w,n}(x)|\cdot\mu(I_w^n).
$$
Denote the Perron-Frobenius operator acting on the space of densities also by $F_*$. Let $x\in[-1,1]$ and let $y_i\in I_i^1$ be such that $F(y_i)=x$. Then
$$
\aligned
F_*(\rho)(x)&=\sum_{i} \frac{1}{|DF(y_i)|}\cdot \rho(y_i)\\
&=  \sum_{i} \frac{1}{|DF(y_i)|}\int |D\psi_w(y_i)| d\nu\\
&= \lim_{n\to\infty} \sum_{i}  \sum_{I_w^n\in \Cal{I}^n}  \frac{1}{|DF(y_i)|}\cdot |D\psi_w^n(y_i)|\cdot \mu(I_w^n)\\
&=\lim_{n\to\infty}  \sum_{i}  \sum_{I_w^n\in \Cal{I}^n} |D\psi_{iw,n+1}(x)|\cdot \frac{|I_{iw}^{n+1}|}{|I_w^n|}\cdot \mu(I_w^n),\\
\endaligned
$$
where $iw$ is the word obtained by concatenating $w$ after $i$. Now we will use Lemma 2.5:
$$
\aligned
F_*(\rho)(x)
&=\lim_{n\to\infty}  \sum_{i}\sum_{I_w^n\in \Cal{I}^n} |D\psi_{iw, n+1}(x)|\cdot \mu(I_{iw}^{n+1})\cdot \frac{\mu(I_w^n)}{\mu(I_{iw}^{n+1})}
                 \cdot \frac{|I_{iw}^{n+1}|}{|I_w^n|}\\
&=\lim_{n\to\infty}  \sum_{i}\sum_{I_w^n\in\Cal{I}^n} |D\psi_{iw,n+1}(x)|\cdot \mu(I_{iw}^{n+1})\\
&=\lim_{n\to\infty} \sum_{I_v^{n+1}\in\Cal{I}^{n+1}} |D\psi_{v,n+1}(x)|\cdot \mu(I_{v}^{n+1})\\
&=\rho(x).\\
\endaligned
$$
\hfill\hfill\qed $\,\,$ (Lemma 2.8)
\enddemo

\demo{Proof of Proposition 2.2} Let $H$ be a multiplier preserving conjugacy between the Markov-maps $F$ and $G$. Lemma 2.4 states that each branch
$F^n:I_w^n\to[-1,1]$ and $G^n:H(I_w^n)\to[-1,1]$ has uniformly bounded distortion, say bounded by $K$. Lemma 2.4(3) implies that each branch $F^n:I_w^n\to[-1,1]$ has an unique fixed point $p_w^n$. Now we can estimate the length of $I_w^n$:
$$
\frac{1}{K}\le \frac{2/|I_w^n|}{|DF^n(p_w^n)|}\le K.
$$
And 
$$
\frac{1}{K}\le \frac{2/|H(I_w^n)|}{|DG^n(H(p_w^n))|}{}\le K.
$$
Because $DG^n(H(p_w^n))=DF^n(p_w^n)$ we get
$$
\frac{1}{K}\le \frac{|H(I_w^n)|}{|I_w^n|}\le K,
$$
for each branch $I_w^n\in\Cal{I}^n$.

\flushpar
To prove that $H$ is Lipschitz choose $x,y\in [-1,1]$. Let 
$$
\Cal{I}^n_{x,y}=\{I\in\Cal{I}^n| I\subset [x,y]\}.
$$
By using that $|I_w^n|, |H(I_w^n)|\to 0$ for $n\to \infty$ we get
$$
|H(y)-H(x)|=\lim_{n\to\infty} \sum_{I\in \Cal{I}^n_{x,y}} |H(I)|\le
            \lim_{n\to\infty} K\cdot \sum_{I\in \Cal{I}^n_{x,y}} |I|=K|y-x|.
$$
\hfill\hfill\qed $\,\,$ (Proposition 2.2)
\enddemo

\demo{Proof of Proposition 2.3} Let $H:[-1,1]\to [-1,1]$ be the $C^r$ diffeomorphism defined by
$$
H(x)=\int_{-1}^x \rho(t)dt.
$$
Indeed, $H$ is a diffeomorphism because $\rho$ is bounded and away from zero (see Lemma 2.5). By construction we have $H_*(\mu)=\frac12 \lambda$. In particular, the $C^r$ Markov-map $F_0=H^{-1}\circ F\circ H$ preserves the Lebesgue measure.
\hfill\hfill\qed $\,\,$ (Proposition 2.3)
\enddemo

\proclaim{Corollary 2.9} Every conjugacy between two $C^r$ Markov-maps
 which preserves multipliers is a $C^r$ diffeomorphism. 
\endproclaim

\demo{Proof} There exist $C^r$ Markov-maps $F_0$ and $G_0$ which preserve the Lebesgue measure and are $C^r$ conjugated to respectively $F$ and $G$. The homeomorphism $H$ which preserves the multipliers of $F$ and $G$ induces an exponent preserving conjugacy $H_0$ between $F_0$ and $G_0$. This conjugacy is absolutely continuous, by Proposition 2.2. The ergodicity of both maps $F_0$ and $G_0$, see Lemma 2.4, imply that the conjugacy $H_0$ has to map the invariant absolutely continuous measure of $F_0$ to the absolutely continuous invariant measure of $G_0$. However, both measure are the Lebesgue measure. Consequently, $H_0$ is the identity. 

\hfill\hfill\qed $\,\,$ (Corollary 2.9)
\enddemo

\flushpar
The proof of Corollary 2.9 summarizes the arguments in [SS].

\bigskip
\centerline{\bf 3. Proof of the Theorem}
\bigskip

\bigskip

\flushpar
In this section we will prove the Theorem. Fix two   non-renormalizable unimodal map $f,g\in \Cal{U}$ which does not have periodic attractors neither a Cantor attractor. Assume that they have the same exponents. In particular they are conjugated, say by the homeomorphism $h$. In this case the attractor of $f$ and $g$ will be an interval, respectively $A_f=[f^2(0),f(0)]$ and $A_g=[g^2(0),g(0)]$.

\proclaim{Lemma 3.1} For every interval $T\subset A_f$ there exist $n\ge 0$ and
an interval $M\subset T$ such that $f^n:M\to A_f$ is monotone and onto.
\endproclaim

\demo{Proof} The map $f|A_f$ is not renormalizable. Hence it is topologically
conjugated to a continuous expanding piecewise affine map $f_0:[c_2,c_1]\to [c_2,c_1]$ with critical point $0$ and $c_2=f^2(0)$, $c_1=f(0)$. In [GMT] 
(Lemma 2.6) it has been shown that for each interval $T$ there exists an interval $M\subset T$ and $m\ge 0$ such that $f^m:M\to [c_2,c_1]$ is monotone and $f^m(M)=[c_2,0]$ or $f^m(M)=[0,c_1]$. 
The Lemma follows immediately because $[f_0(c_2), f_0(0)]$ 
contains the expanding fixed point of $f_0\in [0,c_1]$.
 \hfill\hfill\qed $\,\,$ (Lemma 3.1)
\enddemo

\flushpar
Take $x\in \text{int}(A_f)$, the interior of $A_f$. We will describe an {\it induced map} on a neighborhood of $x$ and obtain a topologically defined $C^2$ Markov-map. The periodic orbits of $f$ are dense in $A_f$. This allows us to choose an interval $U=[a,b]$ with $x\in \text{int}(U)$ and $a,b$ are periodic points whose orbit do not intersect $\text{int}(U)$. 

\flushpar
An interval $I\subset A_f$ is called {\it good} if there exist $n_I> 0$ and an interval $T\supset I$ such that $f^{n_I}:(I,T)\to (U,A_f)$ is monotone and onto. Clearly, the number $n_I$ is uniquely defined. The interval $T$ is called the {\it extension} of $I$. 
As a consequence of Lemma 3.1 we get that in every open set there are good intervals.

\flushpar
The fact that the boundary points of $U$ do not return into $\text{int}(U)$ imply that good intervals are nested: if $I_1$ and $I_2$ are good intervals then
$$
\text{int}(I_1)\cap \text{int}(I_2)=\emptyset \text{ or } I_1\subset I_2 \text{ or } I_2\subset I_1.
$$
We say that a good interval $I$ is of {\it first generation} if there does not exist a good interval $J$ which strictly contains $I$. The density of good intervals imply that also the first generation good intervals are dense.

\flushpar
Let $\Cal{I}$ be the collection of first generation good intervals which are part of $U$. Observe that there are no good intervals which intersect the boundary of $U$ in their interior. In general $\Cal{I}$ will be a countable collection and $\overline{\bigcup \Cal{I}}=U$. Define $M:\bigcup \Cal{I}\to U$ by $M:I=f^{n_I}$, $I\in\Cal{I}$.

\proclaim{Lemma 3.2} $|U\setminus \bigcup \Cal{I}|=0$.
\endproclaim

\demo{Proof} We will use some notation from [Ma]. Let $n\ge 0$ and $x\in A_f$. The maximal interval on which $f^n$ is monotone is denoted by $T^n(x)$. The components of $f^n(T_n(x))-\{f^n(x)\}$ are denoted by $L_n(x)$ and $R_n(x)$. Define $r_n:A_f\to \Bbb{R}$ by
$$
r_n(x)=\min\{|L_n(x)|,|R_n(x)|\}.
$$
Furthermore define $r:A_f\to \Bbb{R}$ by
$$
r(x)=\limsup_{n\to \infty} r_n(x).
$$
It has been proved in [Ma] that there exists some $r_f>0$ such that
$$
r(x)=r_f
$$
for almost every $x\in A_f$. Here we used that $f$ does not have a Cantor attractor.

\flushpar
Now we will describe an extra condition on the periodic boundary points $a$ and $b$ of $U$. 
The two periodic orbits define a partition of $A_f$:
$\Cal{P}=\{P_i|i=1,2,\cdots,s\}$ is formed by the intervals between the points of the two periodic orbits. Assume that the maximal length of the intervals in this partition is smaller that $\frac{1}{10}r_f$. But  also that the intervals in the partition $h(\Cal{P})$ have length smaller than $\frac{1}{10}r_g$.

\flushpar
In each interval $P_i$, $i=1,2,\cdots s$ there exists a good interval $G_i\subset P_i$.  The partition $\Cal{P}$ is finite hence there exists an $\epsilon>0$ such that
$$
\frac{|G_i|}{|P_i|}\ge \epsilon, i=1,2,\cdots, s.
$$

\bigskip

\flushpar
We will prove the Lemma by contradiction. Assume that there is a density point $y\in U$ of $|U\setminus \bigcup \Cal{I}|=0$ with $r_f(y)=r$. The point $y$ being a density point implies that for good intervals $I$ close to $y$ 
$$
\lim_{I\to y} \frac{|I|}{\text{dist}(y,I)}=0.
$$
We will construct big good intervals arbitrarily close to $y$ and obtain a contradiction to the above statement.

\flushpar
Because $r_f(y)=r$ there exists a sequence $k_n\to\infty$ such that $r_{k_n}(y)\ge \frac12 r_f$. For each $n\ge 0$ there exists $P_{i_n}\in\Cal{P}$ such that $f^{k_n}(y)\in P_{i_n}$. Moreover, there exists an interval $M_n\subset T_{k_n}(y)$ with $y\in M_n$ and $f^{k_n}:M_n\to P_{i_n}$ monotone and onto. Because $|P_{i_n}|\le \frac{1}{10} r_f$ and $r_{k_n}(y)\ge \frac12 r_f$ we get from the Koebe-Lemma that  $f^{k_n}:M_n\to P_{i_n}$ has uniformly
bounded distortion. 

\flushpar
Let $G'_n\subset M_n$ be such that $f^{k_n}(G'_n)=G_{i_n}$. The uniformly bounded distortion of $F^{k_n}|M_n$ implies
$$
\frac{|G'_n|}{\text{dist}(y,|G'_n|)}\ge \epsilon'>0
$$
for some uniform $\epsilon'>0$. Because $|T_{k_n}(y)|\to 0$ we get a contradiction.

\hfill\hfill\qed $\,\,$ (Lemma 3.2)
\enddemo

\flushpar
Let $V\supset U$ be the maximal interval which contains only the points $a$ and $b$ of the periodic orbits of $a$ and $b$. Clearly, for every $I\in\Cal{I}$ there exists $E_I\supset I$ with $E_I\subset U$ such that
$$
f^{n_I}:E_I\to V
$$
is monotone onto.

\proclaim{Lemma 3.3} There exists $K>0$ such that
$$
|[(M|I)^{-1}]|_2\le K
$$
for every $I\in \Cal{I}$.
\endproclaim

\demo{Proof} We have to show that for some $K>0$  
$$
|\eta_{[(M|I)^{-1}]}|_0\le K, I\in \Cal{I}.
$$
Observe, that all the $[(M|I)]$ have, by construction, a monotone extension with negative Schwarzian derivative up to a fixed interval strictly bigger than $[-1,1]$. The Koebe-Lemma implies that for some $K>0$
$$
\ln(\frac{|DM(x)|}{|DM(y)|})\le K\cdot |x-y|,
$$
for any $x,y\in [-1,1]$. Hence,
$$
|\ln(\frac{[(M|I)^{-1}](x)}{[(M|I)^{-1}](y)})|\le K\cdot |x-y|,
$$
for any $x,y\in [-1,1]$. Consequently, 
$$
|\int_x^y \eta_{|[(M|I)^{-1}]}| =|\ln(\frac{[(M|I)^{-1}](x)}{[(M|I)^{-1}](y)})|\le K\cdot |x-y|,
$$
for any $x,y\in [-1,1]$. This implies that $|\eta_{[(M|I)^{-1}]}|_0\le K$.
\hfill\hfill\qed $\,\,$ (Lemma 3.3)
\enddemo

\flushpar
The previous Lemmas allow us to construct, by rescaling range and domain of $M$, a $C^2$ Markov-map. 
The definition is essentially topological: we get a corresponding $C^2$ Markov-map on the domain $h(U)$ which is induced by $g$. The restriction $h|U$ of the multipliers preserving conjugacy $h$ between $f$ and $g$ becomes after rescaling a multiplier preserving conjugacy between the Markov-maps. Corollary 2.9 states that $h$ is a $C^2$ diffeomorphism. We proved that the conjugacy between $f$ and $g$ is $C^2$ when restricted to $\text{int}(A_f)$.

\flushpar
Observe that $f(f^2(0))\in \text{int}(A_f)$. Now pull back the conjugation around $f(f^2(0))$ two time to show that $h$ is also smooth in $f(0)$ and $f^2(0)$. 
We showed that $f$ and $g$ are smoothly conjugated on a neighborhood of their attractors.

\bigskip

\flushpar
Observe that in general these smooth conjugacy between the attractors can not be smoothly extended to the whole domain of the maps.

\bigskip
\centerline{\bf  Appendix}
\bigskip

\bigskip

\flushpar
The Schwarzian derivative of a $C^3$ map $f:[-1,1]\to [-1,1]$ is
$$
Sf(x)=\frac{D^3f(x)}{Df(x)}-\frac32\cdot(\frac{D^2f(x)}{Df(x)})^2.
$$
where $D^if(x)$ stands for the $i^{th}$ derivative of $f$ in $x\in [-1,1]$.

\proclaim{Koebe-Distortion-Lemma} There exist continuous functions $K:(0,\infty) \to (0,\infty)$ and $\delta:(0,\infty) \to (0,\infty)$ with the following properties. Let $f^n:T\to f^n(T)$ be monotone and $Sf(x)\le 0$
for all $x\in [-1,1]$. If $I\subset T$ is an interval such that both component
$L,R\subset T\setminus I$ satisfy
$$
\frac{|f^n(L)|}{|f^n(I)|},\frac{|f^n(R)|}{|f^n(I)|}\ge \tau
$$
then $f^n|I$ has bounded {\it distortion}:
$$
\frac{|Df^n(x)|}{|Df^n(y)|}\le 1+K(\tau)\frac{|x-y|}{|I|},
$$
for all $x,y\in I$ and
$$
\frac{|L|}{|I|},\frac{|R|}{|I|}\ge \delta(\tau).
$$
Moreover, $\lim_{\tau\to\infty} K(\tau)=1$ and $\lim_{\tau\to\infty} \delta(\tau)=\infty$.
\endproclaim

\flushpar
A detailed discussion of the Koebe-Lemma can be found in [MS].

\bigskip

\flushpar
The Lebesgue measure on $[-1,1]$ is denoted by $|.|$. A smooth map on the interval is called {\it ergodic} if the only measurable forward invariant sets have Lebesgue measure zero of full measure.

\bigskip
\centerline{\bf References}
\bigskip

\parindent=40pt
\item{[BKNS]} H.Bruin, G.Keller, T.Nowicki, S.van Strien, {\it Wild
               Attractors Exist}, Ann. of Math. {\bf 143} (1996),  

\item{[BL]} A.M.Blokh, M.Lyubich, {\it Attractors of Maps on the interval},
            Func.Anal. and Appl. {\bf 21} (1987), 32-46.

\item{[GMT]} R.Galeeva, M.Martens, C.Tresser, {\it Inducing, Slopes, and        
             Conjugacy Classes}, Israel Jour.Math. {\bf 99} (1997), 123-147.

\item{[L1]}  M.Lyubich, {\it Teichm\"uller Space of Fibonacci Maps}, 
             Stony Brook preprint 1993/12.

\item{[L2]} M.Lyubich, {\it Combinatorics, Geometry and Attractors of
            Quasi-Quadratic Maps}, Ann. of Math {\bf 140} (1994), 347-404.

\item{[LM]} M.Lyubich, J.W.Milnor, {\it The Fibonacci Unimodal Map},
            J.A.M.S. {\bf 6} (1993), 425-457.

\item{[Ma]} M.Martens, {\it Distortion Results and Invariant Cantor Sets of
           Unimodal Maps}, Erg.Th. \& Dyn.Sys. {\bf 14} (1994), 331-349.

\item{[Mi]} J.W.Milnor, {\it On the Concept of Attractor}, Commun. Math. Phys.
             {\bf 99} (1985), 177-195.

\item{[MS]} W.de Melo and S.van Strien, {\it One-dimensional Dynamics},
            Springer-Verlag, 1993.

\item{[SS]} M.Shub, D.Sullivan {\it  Expanding Endomorphism of the Circle  
            revisited}, Erg.Th. \& Dyn.Sys. {\bf 5} (1985), 285-289.

\bye